\documentclass[11pt]{article}
\usepackage{amsmath,epsfig,latexsym,amssymb,cite,graphicx,a4}
\newlength{\indentedwidth}
\newdimen\mathindent
\indentedwidth=\mathindent

\usepackage{epsf,graphics}
\DeclareMathAlphabet{\mathpzc}{OT1}{pzc}{m}{it}  
\begin{document}
\vskip 0.5cm
\begin{center}
{\Large \bf Universal Baxterization for $\mathbb{Z}$-graded Hopf algebras}
\end{center}
\vskip 0.8cm
\centerline{K.A. Dancer%
\footnote{\tt dancer@maths.uq.edu.au}, 
P.E. Finch%
\footnote{\tt pfinch@maths.uq.edu.au}
and P.S. Isaac%
\footnote{\tt psi@maths.uq.edu.au}}
\vskip 0.9cm
\centerline{\sl\small 
 Centre for Mathematical Physics, School of Physical Sciences, }
\centerline{\sl\small The University of Queensland, Brisbane 4072,
 Australia.} 
\vskip 0.9cm
\begin{abstract}
We present a method for Baxterizing solutions of the constant Yang-Baxter
equation associated with $\mathbb{Z}$-graded Hopf algebras. To demonstrate the
approach, we provide examples for the Taft algebras and the quantum group
$U_q\left[sl(2)\right]$. 
\end{abstract}

\setcounter{footnote}{0}
\def\thefootnote{\fnsymbol{footnote}}

\newtheorem{definition}{Definition}[section]
\newtheorem{lemma}{Lemma}[section]
\newtheorem{proposition}{Proposition}[section]
\newtheorem{theorem}{Theorem}[section]
\newtheorem{example}{Example}[section]

\def\C{\mathbb{C}}
\def\Z{\mathbb{Z}}
\def\a{\alpha}
\def\b{\beta}
\def\e{\varepsilon}
\def\d{\delta}
\def\g{\gamma}
\def\l{\lambda}
\def\k{\kappa}
\def\End{\mbox{End\ }}

\section{Introduction}

The word ``Baxterization'' was originally coined by V.F.R. Jones
\cite{Jones,Jones2} to refer
to the insertion of a parameter into a solution of the constant Yang-Baxter equation so that
it becomes a solution of the parameter dependent Yang-Baxter equation. This is
done in such a way that the resultant parametric solution reduces to the
original constant one in some suitable limit.

There exist well-studied methods of Baxterization, especially those associated with 
quantum groups.
Both universal (i.e. representation independent) \cite{Jones2,KRS,KulishResh,ZhangGould} and
representation dependent \cite{Jones,Jones2,Jimbo,Jimbo3,ZGB,DGZ,ChangGeXue,ZhangKauffmanGe,WellyMartins,DIL} 
approaches have been developed. 

In this paper, we introduce a new method of Baxterizing universal R-matrices
arising from $\mathbb{Z}$-graded associative algebras. In particular
we focus on $\mathbb{Z}$-graded Hopf algebras. The prime examples with which we demonstrate
our results are the Taft algebras \cite{Taft,Chen}.

The (multiplicative) parameter dependent Yang-Baxter equation (YBE) is 
$$
R_{12}(x)R_{13}(xy)R_{23}(y) = R_{23}(y)R_{13}(xy)R_{12}(x).
$$
Here $R$, known as an $R$-matrix, is an operator on $V \otimes V$ for some vector space $V$.  
We use the standard notation that $R_{13}\in \End(V\otimes V\otimes V)$ 
represents $R$ operating on the $1$st and $3$rd components
of $V\otimes V\otimes V$, and similarly for $R_{12}, R_{23}$.
This equation has a variety of applications, 
particularly in exactly solvable models in statistical mechanics \cite{Bax} and
quantum field theory \cite{Faddeev}. 
Consequently, it is always of interest to develop new methods of solving this
equation. Certainly there already exists a body of elegant works dedicated to
solving this equation, some noteworthy articles being
\cite{Drinfeld,Drinfeld2,Jimbo}. For a good overview of the parameter
dependent Yang-Baxter equation and its
solutions, see for example \cite{Jimbo2} or \cite{Lambe}.

By contrast, the constant Yang-Baxter equation 
$$
R_{12}R_{13}R_{23} = R_{23}R_{13}R_{12}
$$
has no parameter dependence and hence is easier to solve. Solutions are known in
many different contexts, most significantly Drinfeld's universal solution
arising from the quantum double construction for Hopf algebras \cite{Drinfeld}. 
From a solution of the parameter dependent YBE, one can easily obtain a solution 
of the constant YBE (by taking some suitable limit), but the converse is not
true. As mentioned above, there are well established Baxterization techniques
for quantum groups, however these methods do not extend to Hopf algebras in
general.

In this paper we present a straightforward method of obtaining a universal parameter
dependent solution from a constant solution in the context of
$\mathbb{Z}$-graded Hopf algebras.
To demonstrate the method, we provide specific examples for the finite dimensional Taft algebras and
the quantum group $U_q\left[sl(2)\right]$. 

\section{Universal Baxterization}

\begin{definition} 
Let $H$ be an associative algebra with unit with multiplication $m$.  Let $A$ be a subalgebra of $H$.
If we can find $\{ A^p| p \in \mathbb{Z}\}$ such that

\begin{itemize} 
\item[(i)] $A = \bigoplus_{p} A^p$
\end{itemize}
and
\begin{itemize}
\item[(ii)] $m: A^p \otimes A^q \rightarrow A^{p+q},$ 
\end{itemize}
then we say $A$ is $\mathbb{Z}$-graded, and call $A = \bigoplus_p A^p$ the $\mathbb{Z}$-grading of $A$.  If there exists some $p \neq 0$ such 
that $A^p \neq \{ 0\}$ we say the $\mathbb{Z}$-grading is nontrivial.

\end{definition}

\begin{proposition} \label{Baxterization}
Let $H$ be an associative algebra with unit.  Suppose $H$ has subalgebras $A,B$ with $\mathbb{Z}$-gradings 
$A = \bigoplus_p A^p$ and $B= \bigoplus_q B^q$ respectively.  If $H$ contains a
solution of the constant Yang--Baxter equation of the form

$$R = \sum_{i,\a} a_\a^i \otimes b_\a^i$$
where $a_\a^i \in A^i$ and $b_\a^i \in B^i$, then 

$$R(\mu) = \sum_i \mu^i \sum_\a a_\a^i \otimes b_\a^i$$
is a solution of the multiplicative parametric Yang--Baxter equation. \\
{\bf Proof:} 
It is given that $R = \sum_{i,\a} a_\a^i \otimes b_\a^i$ satisfies the constant Yang-Baxter equation

$$R_{12}R_{13} R_{23} = R_{23} R_{13} R_{12}.$$
\noindent
Substituting in, this is equivalent to stating

$$\sum_{i,j,k,\a,\b,\g} a^i_\a a^j_\b \otimes b^i_\a a^k_\g \otimes b^j_\b b^k_\g = 
\sum_{p,q,r,\d,\e,\k} a^q_\e a^r_\k \otimes a^p_\d b^r_\k \otimes b^p_\d b^q_\e.$$
\noindent
In particular, we can equate the entries belonging to $A^s \otimes H \otimes B^q$, giving

$$ \sum_{j,\a,\b,\g} a^{s-j}_\a a^j_\b \otimes b^{s-j}_\a a^{t-j}_\g \otimes b^j_\b b^{t-j}_\g = 
\sum_{q,\d,\e,\k} a^q_\e a^{s-q}_\k \otimes a^{t-q}_\d b^{s-q}_\k \otimes b^{t-q}_\d b^q_\e.$$
\noindent
Now we substitute the parametrized $R$-matrix $R(\mu)$ into the parametric Yang-Baxter equation:

\begin{align*}
R_{12}(\mu) R_{13}(\mu \nu) R_{23}(\nu) 
&= \sum_{i,j,k,\a,\b,\g} \mu^{i+j}\nu^{j+k} \, a^i_\a a^j_\b \otimes b^i_\a a^k_\g \otimes b^j_\b b^k_\g \\
&= \sum_{s,t} \mu^s \nu^t \sum_{j,\a\b\g} a^{s-j}_\a a^j_\b \otimes b^{s-j}_\a a^{t-j}_\g \otimes b^j_\b b^{t-j}_\g \\
&= \sum_{s,t} \mu^s \nu^t \sum_{q,\d,\e,\k} a^q_\e a^{s-q}_\k \otimes a^{t-q}_\d b^{s-q}_\k \otimes b^{t-q}_\d b^q_\e \\
&= \sum_{p,q,r,\d,\e,\k} \mu^{q+r} \nu^{p+q} \, a^q_\e a^r_\k \otimes a^p_\d b^r_\k \otimes b^p_\d b^q_\e \\
&= R_{23}(\nu) R_{13}(\mu \nu) R_{12}(\mu) 
\end{align*}

\noindent as required.

\end{proposition}

It is possible to generalize the result of Proposition \ref{Baxterization} to
include $\mathbb{Z}^n$-graded algebras. The result is that if
$R$ is an element of $\bigoplus_{p\in\mathbb{Z}^n}A^p\otimes B^p$, then a
universal Baxterization exists. Explicitly, let
$$
R = \sum_{i\in\mathbb{Z},\ p\in\mathbb{Z}^n} a_i^p\otimes b_i^p,
$$
where $a_i^p\in A^p,$ $b_i^p\in B^p$. The Baxterized solution will then be
$$
R(\mu) =
\sum_{p\in\mathbb{Z}^n}\mu^{\tau(p)}\sum_{i\in\mathbb{Z}}a_i^p\otimes b_i^p,
$$
where $\tau:\mathbb{Z}^n\rightarrow \mathbb{Z}$ (or some other appropriate codomain) is a group homomorphism under
addition. The proof of this result is essentially the same as the proof of
Proposition \ref{Baxterization}. We will not, however, make use of this
generalization in the current article.

One algebraic structure where a nontrivial $\mathbb{Z}$-grading may arise is the Drinfeld double of a Hopf algebra.  
To understand the Drinfeld double, we first introduce the dual of a finite Hopf algebra $H$, which we denote $H^*$.

The vector space underlying $H^*$ is the set of linear maps $f:H \rightarrow \C$.  We choose the bilinear form 

$$\langle f, x \rangle = f(x), \quad \forall x \in H.$$
If $H$ has basis $\{ a_i \}$, then we choose $\{a_i^*\}$ as a basis for $H^*$, where 

$$\langle a_i^*, a_j \rangle = \delta_{ij}.$$ 
The structure of $H^*$ is induced by that of $H$.  Specifically, if $H$ has multiplication $m$, unit $u$, coproduct $\Delta$ 
and counit $\e$ then $H^*$ becomes a Hopf algebra with multiplication 
$m^*$, unit $u^*$, coproduct $\Delta^*$ and counit $\e^*$ defined by:

\begin{alignat*}{2}
&\langle m^*(a_{i}^{*} \otimes a_{j}^{*}), a_{k} \rangle  = \langle a_{i}^{*} \otimes a_{j}^{*}, \Delta(a_{k}) \rangle , & \qquad &
\langle u^*(\textit{k}), a_{i} \rangle   =  k\epsilon(a_{i}),\; \forall k \in \C, \\
&\langle \Delta^*(a_{i}^{*}), a_{j} \otimes a_{k} \rangle   =  \langle a_{i}^{*}, m(a_{j} \otimes a_{k}) \rangle,  & &
\epsilon^*(a_{i}^{*})  =  \langle a_{i}^{*},e \rangle.
\end{alignat*}

The Drinfeld double of a finite Hopf algebra $H$, which we denote $D(H)$, is a quasitriangular Hopf algebra spanned by elements of the form 
$\{gh^*|g \in H, h^* \in H^*\}$.  Details of the algebraic structure and costructure of $D(H)$ can be found in \cite{Drinfeld}.  
Of particular relevance here is the property that $D(H)$ contains a canonical
solution of the Yang--Baxter equation of the form

$$R = \sum_{i} a_i \otimes a_i^*,$$
where $\{a_i\}$ is a basis for $H$.  Here we identify $a_i$ with $a_i \epsilon$ and $(a_i)^*$ with $u(a_i)^*$ where $\epsilon$ and $u$ are the counit and unit of $H$ respectively.

Using this universal $R$-matrix, we have the following result:

\begin{proposition} \label{double}
Let $H$ be a finite-dimensional $\mathbb{Z}$-graded Hopf algebra with nontrivial $\mathbb{Z}$-grading $H = \bigoplus_p A^p$.  
If the coproduct of $H$ satisfies
$$\Delta: A^p \rightarrow \bigoplus_q A^q \otimes A^{p-q}, \qquad \forall p \in \mathbb{Z},$$
then $D(H)$ nontrivially satisfies the conditions for Proposition \ref{Baxterization}.

\

\noindent {\bf Proof:}
Let $A^{p}$ have the basis $\{ a^{p}_{i} \}$ and $B^{p}$ have the basis $\{ (a^{p}_{i})^*|a^p_i \in A^p \}$.  Clearly the dual of $H$ can be written as 
$H^* = \bigoplus_p B^p.$  Moreover an $R$-matrix for $D(H)$ is $\sum_{j} a_j \otimes (a_j)^* \; \in \bigoplus_p A^p \otimes B^p$.  Thus it remains 
only to show $m: B^p \otimes B^q \rightarrow B^{p+q}$ where $m$ represents multiplication within $D(H)$.  

But
\begin{align*} 
\langle m \bigl( (a^p_i)^* \otimes (a^q_j)^* \bigr), a^r_k \rangle 
&= \langle (a^p_i)^* \otimes (a^q_j)^*, \Delta(a^r_k) \rangle \\
&= 0 \hspace{4cm} \mbox{if $r \neq p + q$}.
\end{align*}
Thus $m: B^{p} \otimes B^{q} \rightarrow B^{p+q}$ for all $p, q \in \Z$. Hence $D(H)$ satisfies the conditions for Proposition \ref{Baxterization}.
\end{proposition}

\section{Example: $U_q[sl(2)]$}

The $q$-deformed Lie algebra $U_q[sl(2)]$ has generators $e,f,h$ satisfying

$$[e,f] = \frac{q^h - q^{-h}}{q-q^{-1}}, \quad [h,e] = 2e, \quad [h,f] = -2f,$$
where $q$ is the deformation parameter.  Define $[n]_q$ and $[n]_q!$ as follows:

\begin{align*}
[n]_q &= \frac{q^n - q^{-n}}{q-q^{-1}}, \\
[n]_q! &= [n]_q [n-1]_q ... [1]_q.
\end{align*}
Then $U_q[sl(2)]$ contains the following universal $R$-matrix \cite{CP}:

$$\mathfrak{R} = \sum_{n=0}^\infty \frac{q^{\frac 12 n(n+1)}(1-q^{-2})^n}{[n]_q!} q^{\frac 12 (h \otimes h)} e^n \otimes f^n.$$
Now set $H = \langle h \rangle$ to be the subalgebra generated by $h$, $A = \langle h,e \rangle$ to be the subalgebra generated by 
$e$ and $h$, and $B = \langle h,f \rangle$ to be the subalgebra generated by $f$ and $h$.  Then $A$ has the natural 
$\mathbb{Z}$-grading $A = \bigoplus_{k \in \mathbb{N}} A^k$ where 
$A^k = H e^k, \, k \geq 0$.  Similarly, $B$ has the $\mathbb{Z}$-grading $B = \bigoplus_{k \in \mathbb{N}} B^k$ where $B^k = H f^k, \, k \geq 0$.  
Note that with these $\mathbb{Z}$-gradings $R \in \bigoplus_{k \in \mathbb{N}} A^k \otimes B^k$, so we can apply Proposition \ref{Baxterization}.
 We find that

$$\mathfrak{R}(\mu) = \sum_{n=0}^\infty \mu^n \frac{q^{\frac 12 n(n+1)}(1-q^{-2})^n}{[n]_q!} q^{\frac 12 (h \otimes h)} e^n \otimes f^n$$
is a solution of the mulitplicative parametric Yang-Baxter equation.

Applying the spin-$\frac 12$ representation, which is given by 

\begin{equation*}
e = \begin{pmatrix} 0 & 1 \\ 0 & 0 
\end{pmatrix}, \quad
f = \begin{pmatrix} 0 & 0 \\ 1 & 0
\end{pmatrix}, \quad
h = \begin{pmatrix} 1 & 0 \\ 0 & -1
\end{pmatrix},
\end{equation*}
this becomes

$$R(\mu) = \begin{pmatrix}
q^{\frac 12} & 0 & 0 & 0 \\
0 & q^{-\frac 12} & \mu q^{-\frac 12} (q - q^{-1}) & 0 \\
0 & 0 & q^{-\frac 12} & 0 \\
0 & 0 & 0 & q^{\frac 12}
\end{pmatrix}.$$
Similarly, applying the spin-$1$ representation, which is given by 

\begin{equation*}
e = \sqrt{q + q^{-1}} \begin{pmatrix} 0 & 1 & 0\\ 0 & 0 & 1  \\0 & 0 & 0 
\end{pmatrix}, \quad
f = \sqrt{q + q^{-1}} \begin{pmatrix} 0 & 0 & 0 \\ 1 & 0 & 0 \\ 0 & 1 & 0
\end{pmatrix}, \quad
h = \begin{pmatrix} 2 & 0 & 0 \\ 0 & 0 & 0 \\ 0 & 0 & -2
\end{pmatrix},
\end{equation*}
the parametric $R$-matrix becomes 
{\small
\begin{equation*}
R(\mu) = \begin{pmatrix} 
q^2 & 0 & 0 & 0 & 0 & 0 & 0 & 0 & 0 \\
0 & 1 & 0 & \mu(q^2 - q^{-2}) & 0 & 0 & 0 & 0 & 0 \\
0 & 0 & q^{-2} & 0 & \mu q^{-2}(q^2 - q^{-2}) & 0 & \mu^2 q^{-1}(q-q^{-1})^2(q+q^{-1})& 0 & 0 \\
0 & 0 & 0 & 1 & 0 & 0 & 0 & 0 & 0  \\
0 & 0 & 0 & 0 & 1 & 0 & \mu(q^2 - q^{-2}) & 0 & 0 \\
0 & 0 & 0 & 0 & 0 & 1 & 0 & \mu(q^2 - q^{-2}) & 0 \\
0 & 0 & 0 & 0 & 0 & 0 & q^{-2} & 0 & 0 \\
0 & 0 & 0 & 0 & 0 & 0 & 0 & 1 & 0  \\
0 & 0 & 0 & 0 & 0 & 0 & 0 & 0 & q^2
\end{pmatrix}.
\end{equation*}
}

\section{Example: Taft algebras}
The Taft algebra $T_{N,q}$ \cite{Taft} over a field $\mathbb{F}$ is an $N^2$-dimensional algebra with unit $e$ generated by $\langle a,x|a^N = e, \, x^N = 0, xa = qax \rangle$.  Here $q$ is a primitive $N^{\rm th}$ root of unity in $\mathbb{F}$.  We choose $\{ a^i x^j| 0 \leq i,j < n\}$ as a basis for $T_{N,q}$, and note that multiplication of two basis elements is given by $(a^i x^j)( a^k x^l) = q^{jk} a^{i+k} x^{l+j}$. 

The Taft algebra $T_{N,q}$ becomes a Hopf algebra when endowed with a costructure and antipode defined on the generators $a,x$ by:

\begin{alignat*}{3}
&\Delta(a) = a \otimes a, & & \epsilon(a) = 1,  & & \gamma(a) = a^{-1}, \\ 
&\Delta(x) = x \otimes e + a \otimes x, \qquad && \epsilon(x) = 0, \qquad && \gamma(x) = -a^{-1} x.
\end{alignat*}
Here $\Delta, \, \epsilon$ and $\gamma$ represent the coproduct, counit and antipode respectively.  

The coproduct and counit extend as homomorphisms to all of $T_{N,q}$.  Following the notation of \cite{Chen}, we define 
$(n)_q = 1 + q + ... + q^{n-1}$ and $(n)_q! = (n)_q (n-1)_q ... (1)_1.$  Set 

$${n \choose m}_{q} = \frac{(n)_{q}!}{(m)_{q}!(n-m)_{q}!}.$$

\noindent
Then for all elements $a^i x^j \in T_{N,q}$, we find the coproduct is given by

$$\Delta(a^{i}x^{j}) = \sum_{k=0}^{j} {j \choose k}_{q} a^{j-k+i}x^{k} \otimes a^{i}x^{j-k}.$$

The Drinfeld double $D(T_{n,q})$ contains a universal $R$-matrix given by

$$R = \sum_{i,j=0}^{N-1} a^i x^j \otimes (a^i x^j)^*.$$
\noindent
But $T_{N,q}$ has the $\Z$-grading $T_{N,q} = \bigoplus_p A^p$ where $A^p$ has basis $\{a^i x^p| 0 \leq i < N \}$.  Under this grading, the coproduct 
satisfies $$\Delta: A^p \rightarrow \bigoplus_q A^q \otimes A^{p-q} \quad \forall p \in \Z.$$

Thus we note from Propositions \ref{Baxterization} and \ref{double} that the
Drinfeld double $D(T_{N,q})$ contains an algebraic solution of 
the parametric Yang-Baxter equation given by:

$$R(\mu) = \sum_{i,j=0}^{N-1} \mu^j a^i x^j \otimes (a^i x^j)^*.$$

This can in turn give rise to several matrix solutions of the parametric Yang-Baxter equation, as the representation theory of the Taft algebras has 
been developed by Chen \cite{Chen}.  Explicitly, the $N^2$ irreducible representations of $T_{N,q}$ are given by

$$ \pi_{n,l}(a^{i}x^{j}) = \sum_{k = 1}^{n-j} q^{(k-l-n)i} \frac{(k+j-1)_{q}!}{(k-1)_{q}!} \Pi_{p=0}^{j-1}(1-q^{p+k-n}) e_{k,k+j}$$
and
$$\pi_{n,l}((a^{i+l-1}x^{j})^{*}) = \left\{
\begin{array}{lcl}
	 \frac{e_{i+j,i}}{(j)_{q}!},  &   & 1 \leq i \leq n - j \mbox{  mod(N)},\\
	0 &\quad & \mbox{otherwise},
\end{array} \right.$$
where $1 \leq n,l \leq N.$  Here $e_{i,j}$ is the $n \times n$-dimensional elementary matrix whose only non-zero entry is a 1 in the $(i,j)$ position.

There are also $N$-dimensional indecomposable representations of $T_{N,q}$, which can be found in \cite{Chen}. They are given by:

$$\pi_\a (a^{i}x^{j}) = \alpha q^{-i(j+l)} \frac{(N-2)_{q}!}{(N-j-1)_{q}!} \prod_{p=1}^{j-1} (1 - q^{-p}) e_{N+1-j,1} + \sum_{k=2}^{N-j}q^{i(k-1-l)}\frac{(k+j-1)_{q}!}{(k-2)_{q}!} \prod_{p=0}^{j-1} (1 - q^{k+p}) e_{k,k+j}$$
and, as before,

$$\pi_\a((a^{i+l-1}x^{j})^{*}) = \left\{
\begin{array}{lcl}
	 \frac{e_{i+j,i}}{(j)_{q}!},  &   & 1 \leq i \leq n - j \mbox{  mod(N)},\\
	0 &\quad & \mbox{otherwise}.
\end{array} \right.$$

When the representation \cite{Chen} arising from $\tilde{V}_{3,l} \otimes \tilde{V}_{3,l}$ is applied to the universal R-matrix of $D(T_{N,q})$, $N\geq 3$, it gives the Baxterized $R$-matrix
{\small $$
	R(\mu) =
	\left(
	\begin{array}{ccc ccc ccc}
		1 & 0 & 0 & 0 & 0 & 0 & 0 & 0 & 0 \\
		0 & q^{-l-2} & 0 & (1-q^{-2}) \mu & 0 & 0 & 0 & 0 & 0 \\
		0 & 0 & q^{-2(l+2)} & 0 & q^{-l-4}(q^{2}-1)\mu & 0 & (1-q^{-1})(1-q^{-2})\mu^{2} & 0 & 0 \\
		0 & 0 & 0 & q^{l} & 0 & 0 & 0 & 0 & 0 \\
		0 & 0 & 0 & 0 & q^{-1} & 0 & q^{l+1}(1-q^{-2})\mu & 0 & 0 \\
		0 & 0 & 0 & 0 & 0 & q^{-l-2} & 0 & (1-q^{-2})\mu & 0 \\
		0 & 0 & 0 & 0 & 0 & 0 & q^{2l} & 0 & 0 \\
		0 & 0 & 0 & 0 & 0 & 0 & 0 & q^{l} & 0 \\
		0 & 0 & 0 & 0 & 0 & 0 & 0 & 0 & 1 \\
	\end{array}
	\right)
$$}
It is worth noting that for $l=N-1$, this matrix is a special case of the $9\times 9$ matrix associated with $U_q[sl(2)]$ given earlier.

Although the examples we have given are all upper triangular, this is
not true in general. For a resultant matrix which is not triangular, we require a
$\mathbb{Z}$-grading which is not an $\mathbb{N}$-grading. For example, 
the quantum double of $U_q[sl(2)]$ contains such a grading.

\end{document}